\documentclass[preprint,11pt]{elsarticle}

\usepackage{amssymb}
\usepackage{amsmath}
\usepackage{amsfonts}
\usepackage{amssymb}
\usepackage{indentfirst,latexsym,bm}
\usepackage{amsthm}
\usepackage{color,xcolor}
\usepackage[all]{xy}
\input{mathrsfs.sty}
\newtheorem{theorem}{Theorem}[section]

\newtheorem{proposition}[theorem]{Proposition}

\theoremstyle{definition}

\theoremstyle{definition}

\theoremstyle{remark}

\theoremstyle{question}

\numberwithin{equation}{section}

\journal{XXX}

\begin{document}

\begin{frontmatter}



\title{Some simplified formulas for the matched projection of an idempotent}
\author{Qingxiang Xu}
\ead{qingxiang\_xu@126.com}
\address{Department of Mathematics, Shanghai Normal University, Shanghai 200234, PR China}

\begin{abstract} Let $\mathcal{L}(H)$ be the set of all adjointable operators on a Hilbert $C^*$-module $H$. For each $T\in\mathcal{L}(H)$, $T^*$ denotes its adjoint operator, and $|T|$ is the positive square root of $T^*T$. We establish simplified  formulas for the matched projection $m(Q)$ of an idempotent $Q\in\mathcal{L}(H)$ as
\begin{align*}m(Q)=&\frac{I+|Q^*|-|I-Q^*|}{2}=\frac{I+|Q|-|I-Q|}{2}\\
=&\frac{2+|Q+Q^*|-|2-(Q+Q^*)|}{4},
\end{align*}
where $I$ is the identity operator on $H$. These explicit expressions facilitate the straightforward derivation of several known properties of $m(Q)$.
\end{abstract}


\end{frontmatter}



\section{Introduction}

Throughout this paper, $H$ and $K$ are right Hilbert modules over a $C^*$-algebra.  The set of all adjointable operators from $H$ to $K$   is denoted by $\mathcal{L}(H, K)$, and we abbreviate $\mathcal{L}(H)=\mathcal{L}(H,H)$ when $H=K$. The identity operator on $H$ is denoted by $I_H$, or simply $I$ when no confusion arises. For each $T\in\mathcal{L}(H,K)$, we denote by
$T^*$ its adjoint operator, by $\mathcal{R}(T)$ its range, and by $\mathcal{N}(T)$ its null space. We also use $|T|$ to denote the positive square root of $T^*T$. When $\mathcal{R}(T)$ is closed in $K$, we use the notation $T^\dag$ to denote the Moore-Penrose inverse of $T$.  An operator $Q\in\mathcal{L}(H)$ is called an idempotent if $Q^2=Q$. If, in addition,  $Q^*=Q$, then $Q$ is a projection.

Given an idempotent $Q\in\mathcal{L}(H)$, one may consider the set of operator distances from all projections on $H$ to $Q$. A natural problem is to determine the minimum, maximum, and intermediate values of these distances. To address this problem,
for every idempotent $Q\in\mathcal{L}(H)$, a projection $m(Q)$, called the matched projection of $Q$, was introduced   in \cite{TXF01} as follows:
\begin{equation*}\label{equ:m(q)}
m(Q)=\frac12\big(|Q^*|+Q^*\big)|Q^*|^\dag\big(|Q^*|+I\big)^{-1}\big(|Q^*|+Q\big).
\end{equation*}
Some investigations and applications of the matched projections can be found in \cite{TXF01,TXF02,ZTX01,ZTX02,ZFL}.

In this paper, we show that the expression for the matched projection $m(Q)$ can in fact be simplified to
\begin{equation}\label{equ:simplification}m(Q)=\frac{I+|Q^*|-|I-Q^*|}{2},\end{equation}
which further yields two formulations
\begin{align*}m(Q)=\frac{I+|Q|-|I-Q|}{2}=\frac{2+|Q+Q^*|-|2-(Q+Q^*)|}{4}.
\end{align*}
These explicit expressions enable straightforward derivation of several known properties of $m(Q)$ originally established in \cite{TXF01}.

\section{Derivation of the simplified expression}

Let $Q\in\mathcal{L}(H)$ be an idempotent. We proceed to verify the validity of \eqref{equ:simplification}. First, if $Q\in\mathcal{L}(H)$ is a projection, then
$m(Q)=Q$, $|Q^*|=Q$ and $|(I-Q)^*|=I-Q$. Thus, identity \eqref{equ:simplification} holds trivially.

Now suppose $Q\in\mathcal{L}(H)$ is not a  projection. Let $$H_1=\mathcal{R}(Q),\quad H_2=\mathcal{N}(Q^*).$$
Denote the identity operators on $H_1$ and $H_2$ by $I_1$ and $I_2$, respectively. With respect to the orthogonal decomposition $H=H_1+H_2$, the operator $Q$ admits  the matrix representation
$$Q=\left(
      \begin{array}{cc}
        I_1 & A \\
        0 & 0 \\
      \end{array}
    \right)$$
for some $A\in\mathcal{L}(H_2,H_1)$. This yields
$$|Q^*|=\left(
          \begin{array}{cc}
            B & 0 \\
            0 & 0 \\
          \end{array}
        \right)\quad{and}\quad  (I-Q)(I-Q)^*=\left(
                                     \begin{array}{cc}
                                       AA^* & -A \\
                                       -A^* & I_2 \\
                                     \end{array}
                                   \right),$$
where  $B=(I_1+AA^*)^\frac12\in  \mathcal{L}(H_1)$.
Define $X\in \mathcal{L}(H_2,H_1)$ and $Y\in \mathcal{L}(H_2) $ by
$$ X=A(I_2+A^*A)^{-\frac14}\quad\mbox{and}\quad Y=-(I_2+A^*A)^{-\frac14}.$$
Let $T\in \mathcal{L}(H)$ be the positive operator induced  by
$$T=\left(
      \begin{array}{c}
        X \\
        Y \\
      \end{array}
    \right)\left(
      \begin{array}{c}
        X \\
        Y \\
      \end{array}
    \right)^*=\left(
                \begin{array}{cc}
                  XX^* & XY^* \\
                  YX^* & YY^* \\
                \end{array}
              \right)
    .$$
Notably, $A^*(I_1+AA^*)=(I_2+A^*A) A^*$. This implies that for any continuous function $f$ on $[1, 1+\|A\|^2]$,
$$A^*f(I_1+AA^*)=f(I_2+A^*A)A^*\quad\mbox{and}\quad  f(I_1+AA^*)A=Af(I_2+A^*A). $$
Using these identities, along with the definitions of $X$, $Y$ and $T$, we obtain
$$T=\left(
      \begin{array}{cc}
        B^{-1}AA^* & -B^{-1}A \\
       -A^*B^{-1} & (I_2+A^*A)^{-\frac12} \\
      \end{array}
    \right).$$
It then follows that
\begin{align*}T^2=\left(
                \begin{array}{cc}
                  AA^* & -A \\
                  -A^* & I_2 \\
                \end{array}
              \right).
 \end{align*}
Thus, $T^2=(I-Q)(I-Q)^*$, which implies $T=|(I-Q)^*|$. Consequently,
\begin{align*}\frac{I+|Q^*|-|(I-Q)^*|}{2}=\frac12 \left(\begin{array}{cc}
(B+I_1)B^{-1} & B^{-1}A\\
A^*B^{-1} & A^*\big[B(B+I_1)\big]^{-1}A \\
\end{array}\right)=m(Q),
\end{align*}
where the last equality is established in \cite{TXF01}(see \cite[(3.7)]{TXF01}.

\section{Some applications}
 We derive several properties of the matched projection from its simplified expression \eqref{equ:simplification}.

\begin{proposition}\label{prop:basic properties}{\rm \cite[Theorems 3.6, 3.7 and 3.14]{TXF01}} Let $Q\in\mathcal{L}(H)$ be an idempotent. Then the following statements hold:
\begin{enumerate}
\item[{\rm (i)}] $m(I-Q)=I-m(Q)$.
\item[{\rm (ii)}] If $\pi:\mathcal{L}(H)\to\mathcal{L}(K)$ is a unital $C^*$-algebra morphism, then $m\big(\pi(Q)\big)=\pi\big(m(Q)\big)$.
\item[{\rm (iii)}] $m(Q^*)=m(Q)$ and $Q^*=\big(2m(Q)-I\big)Q\big(2m(Q)-I\big)$.
\item[{\rm (iv)}]  $m(Q)$ and $Q$ as idempotents are homotopy equivalent.
\end{enumerate}
\end{proposition}
\begin{proof}Statements (i) and (ii) follow directly from \eqref{equ:simplification}.

(iii). Define positive operators $T,S\in\mathcal{L}(H)$ by
$$T=|Q^*|+|I-Q|\quad\mbox{and}\quad S=|Q|+|I-Q^*|.$$ A straightforward computation yields
\begin{equation*}T^2=I-Q-Q^*+QQ^*+Q^*Q=S^2,
\end{equation*}
implying $T=S$. Using \eqref{equ:simplification}, we obtain
$$2m(Q^*)-I=|Q|-|I-Q|=|Q^*|-|I-Q^*|=2m(Q)-I,$$ thus $m(Q^*)=m(Q)$.

From \eqref{equ:simplification}, we have
$$Q\big(2m(Q)-I\big)=Q(|Q^*|-|I-Q^*|)=|Q^*|.$$
Replacing $Q$ with $Q^*$ and using the identity $m(Q^*)=m(Q)$ yields
$$Q^*\big(2m(Q)-I\big)=|Q|.$$
Therefore,
$$|Q|\cdot |Q^*|=Q^*\big(2m(Q)-I\big)\big[Q\big(2m(Q)-I\big)\big]^*=Q^*\cdot Q^*=Q^*.$$
This leads to
\begin{align*}\big(2m(Q)-I\big)Q\big(2m(Q)-I\big)=&\big(|Q|-|I-Q|\big)Q\cdot Q\big(|Q^*|-|I-Q^*|\big)\\
=&|Q|\cdot |Q^*|=Q^*.
\end{align*}

(iv). For each $t\in [0,1]$, define $Q_t\in\mathcal{L}(H)$ by
$$Q_t=(1-t)P_{\mathcal{R}(Q)}+tQ,$$
where $P_{\mathcal{R}(Q)}$ is the range projection of $Q$. One can verify that $Q_t$ is an idempotent. Thus,  $Q$ and $P_{\mathcal{R}(Q)}$ as idempotents are homotopy equivalent. Furthermore, since
the matched projection of $P_{\mathcal{R}(Q)}$ is itself, it follows from \eqref{equ:simplification} that $m\big(Q(t)\big)$ is a norm-continuous path of projections in $\mathcal{L}(H)$ which starts at $P_{\mathcal{R}(Q)}$ and ends at $m(Q)$.
Consequently, $P_{\mathcal{R}(Q)}$ and $m(Q)$ are homotopy equivalent. We therefore conclude that $Q$ and $m(Q)$ are homotopy equivalent.
\end{proof}

We provide another new formula for the matched projection as follows.
\begin{theorem}Let $Q\in\mathcal{L}(H)$ be an idempotent. Then
\begin{equation*}m(Q)=\frac{2+|Q+Q^*|-|2-(Q+Q^*)|}{4}.
\end{equation*}
\end{theorem}
\begin{proof}By the proof of Proposition~\ref{prop:basic properties}, we have
\begin{equation*}Q^*=|Q|\cdot |Q^*|\quad\mbox{and}\quad Q=|Q^*|\cdot |Q|.
\end{equation*}
It follows that
\begin{equation*}(|Q|+|Q^*|)^2=Q^*Q+QQ^*+Q^*+Q=(Q+Q^*)^2.
\end{equation*}
Hence,
\begin{equation*}\label{positive T}|Q|+|Q^*|=|Q+Q^*|.\end{equation*}
Replacing $Q$ with $I-Q$ yields
$$|I-Q|+|I-Q^*|=|2-(Q+Q^*)|.$$
Since
\begin{align*}m(Q)=\frac12 \big(m(Q^*)+m(Q)\big)=\frac{2+|Q|+|Q^*|-(|I-Q|+|I-Q^*|)}{4},
\end{align*}
the desired expression for $m(Q)$ follows.
\end{proof}

As a consequence of the preceding theorem, we reestablish some known results from \cite{TXF01}.
\begin{proposition}{\rm \cite[Theorem~3.11, Lemmas~3.12 and 4.5]{TXF01}}  Let $Q\in\mathcal{L}(H)$ be an idempotent. Then the following statements hold:
\begin{enumerate}
\item[{\rm (i)}] $m(Q)Qm(Q)\ge m(Q)$;
\item[{\rm (ii)}]$m(Q)=\frac12(Q+Q^*)\big(m(Q)Qm(Q)\big)^\dag$;
\item[{\rm (iii)}]$\mathcal{R}\big(m(Q)\big)\subseteq \mathcal{R}(Q+Q^*)=\mathcal{R}(|Q|+|Q^*|)$, and $\mathcal{R}\big(m(Q)\big)=\mathcal{R}(Q+Q^*)$ if and only if $Q$ is a projection.
\end{enumerate}
\end{proposition}
\begin{proof}Define $T=Q+Q^*$. Since $Q$ is an idempotent, it is easy to show that $$\sigma(T)\subseteq (-\infty,0]\cup [2,+\infty).$$ Denote by
$\mathfrak{A}$ the commutative $C^*$-subalgebra of $\mathcal{L}(H)$ generated by $T$ and $I$, and let $\Gamma: \mathfrak{A}\to C\big(\sigma(T)\big)$ be the $C^*$-algebra isomorphism induced by the Gelfand transform.
Then $\big(\Gamma(T)\big)(t)=t$ for every $t\in\sigma(T)$, and
\begin{align*}&\Gamma(T_+)(t)=\left\{
                                                                           \begin{array}{ll}
                                                                             0, & \hbox{if $t\in (-\infty,0]\cap\sigma(T)$,} \\
                                                                             t, & \hbox{if $t\in [2,+\infty)\cap \sigma(T)$,}
                                                                           \end{array}
                                                                         \right.\\
&\Gamma\big(m(Q)\big)(t)=\frac{2+|t|-|2-t|}{4}=\left\{
                                                                           \begin{array}{ll}
                                                                             0, & \hbox{if $t\in (-\infty,0]\cap\sigma(T)$,} \\
                                                                             1, & \hbox{if $t\in [2,+\infty)\cap \sigma(T)$.}
                                                                           \end{array}
                                                                         \right.
\end{align*}
It follows that
$$T_+\ge 2m(Q),\quad m(Q)T=Tm(Q)=T_+,\quad m(Q)=T_+ T_+^\dag=T T_+^\dag,$$
where $T_+$ denotes the positive part of $T$, and $T_+^\dag$ is its Moore-Penrose inverse.
Applying Proposition~\ref{prop:basic properties}(iii) yields
$m(Q)Q^*m(Q)=m(Q)Qm(Q)$. Consequently,
\begin{align*}&m(Q)Qm(Q)=\frac12 m(Q)Tm(Q)=\frac12 m(Q)T=\frac12 T_+\ge m(Q),\\
&\frac12(Q+Q^*)\big(m(Q)Qm(Q)\big)^\dag=(Q+Q^*)\big(2m(Q)Qm(Q)\big)^\dag=TT_+^\dag=m(Q).
\end{align*}

Since $\Gamma\big(m(Q)\big)(0)=0$, we see that $m(Q)$ belongs to the $C^*$-subalgebra of $\mathcal{L}(H)$ generated by $T$, which implies $\mathcal{R}\big(m(Q)\big)\subseteq \mathcal{R}(T)$. From the proof of Proposition~\ref{prop:basic properties}, we have
$|Q|+|Q^*|=T\big(2m(Q)-I\big)$, so $\mathcal{R}(T)=\mathcal{R}(|Q|+|Q^*|)$, as $2m(Q)-I$ is a unitary operator.

Now, if $Q$ is a projection, then $m(Q)=Q=Q^*$, and thus $\mathcal{R}\big(m(Q)\big)=\mathcal{R}(T)$. Conversely, suppose $\mathcal{R}\big(m(Q)\big)=\mathcal{R}(T)$, then
$T=m(Q)T=T_+$.
Hence,
$$(I-Q)T_+(I-Q^*)=(I-Q)T(I-Q^*)=0,$$
which gives
$$Q-QQ^*=T(I-Q^*)=T_+(I-Q^*)=0.$$
Therefore, $Q=QQ^*$, implying $Q$ is self-adjoint, so $Q$ is a projection.
\end{proof}


\vspace{5ex}






\begin{thebibliography}{99}





\bibitem{TXF01}X. Tian, Q. Xu and C. Fu, The matched projections of idempotents on Hilbert $C^*$-modules, J. Operator Theory  94 (2025), 219--251.

\bibitem{TXF02}X. Tian, Q. Xu and C. Fu, The Frobenius distances from projections to an idempotent matrix, Linear Algebra Appl. 688 (2024), 21--43.











\bibitem{ZTX01}X. Zhang, X. Tian and Q. Xu, The maximum operator distance from an idempotent to the set of projections, Linear Algebra Appl. 692 (2024), 62--70.


\bibitem{ZTX02}X. Zhang, X. Tian and Q. Xu, Some applications of the matched projections of idempotents, preprint.
arXiv:2404.03433v2



\bibitem{ZFL}C. Zhao, Y. Fang and Y. Li, Some characterizations of the quasi-projection pairs and the matched projections, Banach J. Math. Anal. 19 (2025), no. 2, Paper No. 18.






\end{thebibliography}
\end{document}